\documentclass[twoside,11pt]{article}

\usepackage{amsthm,amsmath,amssymb,epsf}
\usepackage[english]{babel}
\usepackage{epsfig}
\usepackage{amsfonts}
\usepackage{natbib}
\usepackage{latexsym}
\usepackage{bm}
\usepackage[latin1]{inputenc}
\usepackage{amsopn,amstext,amscd}
 \newtheorem{thm}{Theorem}
 \newtheorem{cor}{Corollary}
 
 \newtheorem{prop}{Proposition}
 \newtheorem{defn}{Definition}

\newtheorem{exam}{Example}

\def\d{\mathrm{d}}
\def\ee{\mathrm{e}}

\def\RR{\mathbb{R}}

\def\ss{\mathbf{s}}
\def\hh{\mathbf{h}}
\def\xx{\mathbf{x}}

\def\dd{\mathbf{d}}

\def\expz{\exp(\mathbb{Z}_{+})}
\def\query[#1]{\fbox{\texttt{#1}}}

\def\d{\mathrm{d}}
\def\ee{\mathrm{e}}

\def\RR{\mathbb{R}}
\def\rr{\boldsymbol{r}}

\def\ss{\mathbf{s}}
\def\hh{\mathbf{h}}
\def\xx{\mathbf{x}}

\def\expz{\exp(\mathbb{Z}_{+})}
\def\query[#1]{\fbox{\texttt{#1}}}

\begin{document}

\title{Quasi-arithmetic means of covariance functions with potential applications to space-time data}

\author{{\normalsize
Emilio Porcu$\dagger$\thanks{Corresponding author. {\em Email address:
porcu@mat.uji.es, Fax: +34.964.728429}}, Jorge Mateu$\dagger$ and George
Christakos$\ddagger$}
\smallskip\\
\small $\dagger$ Department of Mathematics, Universitat Jaume I, \\
\small Campus Riu Sec, E-12071 Castellón, Spain \\
\small $\ddagger$ Department of Geography, San Diego State University 5500
\\ \small Campanile Drive
San Diego, CA 92182-4493  \\
}

\date{}
\maketitle

\begin{center}
\includegraphics[width=5cm]{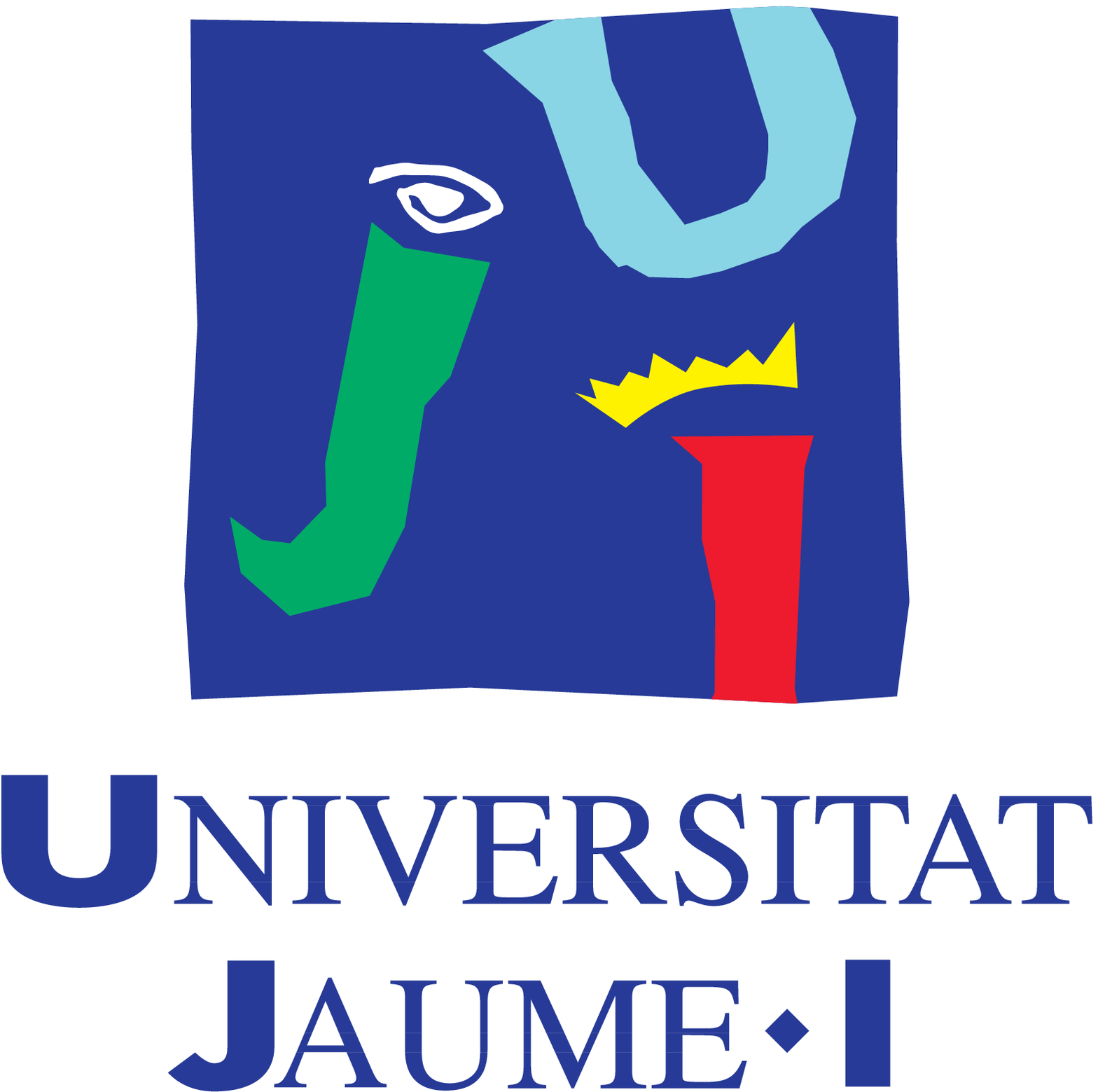} \\
\vspace{2cm} \large {\bf Technical Report n. 107}, {\bf {\em Universitat Jaume I,}} \\
Castellón de la Plana (Spain)
\end{center}
\newpage
\bigskip

\begin{center}
{\bf Abstract}
\end{center}

The theory of quasi-arithmetic means is a powerful tool in the study of
covariance functions across space-time.  In the present study we use
quasi-arithmetic functionals to make inferences about the permissibility
of averages of functions that are not, in general, permissible covariance
functions. This is the case, e.g., of the geometric and harmonic averages,
for which we obtain permissibility criteria. Also, some important
inequalities involving covariance functions and preference relations as
well as algebraic properties can be derived by means of the proposed
approach. In particular, we show that quasi-arithmetic covariances allow
for ordering and preference relations, for a Jensen-type inequality and
for a minimal and maximal element of their class. The general results
shown in this paper are then applied to study of spatial and
spatiotemporal random fields. In particular, we discuss the representation
and smoothness properties of a weakly stationary random field with a
quasi-arithmetic covariance function. Also, we show that the generator of
the quasi-arithmetic means can be used as a link function in order to
build a space-time nonseparable structure starting from the spatial and
temporal margins, a procedure that is technically sound for those working
with copulas. Several examples of new families of stationary covariances
obtainable with this procedure are shown. Finally, we use quasi-arithmetic
functionals to generalise existing results concerning the construction of
nonstationary spatial covariances and discuss the applicability and limits
of this generalisation.

{\bf Subj-class}: Probability; Statistics. \\ {\bf MSC class}:  60G60;
86A32.
\newpage
\section{Introduction}

The importance of quasi-arithmetic means has been well understood at least
since the 1930s, and a number of writers have since contributed to their
characterisation and to the study of their properties. In particular,
Kolmogorov (1930) and Nagumo (1930) derived, independently of each other,
necessary and sufficient conditions for the quasi-arithmeticity of the
mean, that is, for the existence of a continuous strictly monotonic
function $f$ such that, for $x_1, \ldots, x_n$ in some real interval, the
function $(x_1, \ldots, x_n) \mapsto M_n(x_1, \ldots,
x_n)=f^{-1}(\frac{1}{n} \sum_1^{n} f(x_i))$ is a mean.

Using this result, they partially modified the classical Cauchy (1821)
internality and Chisini's (1929) invariance properties. As pointed out by
Marichal (2000), the Kolmogorov reflexive property is equivalent to the
Cauchy internality, and both are accepted by statisticians as requisites
for means.

Early works on the concept of mean include Bonferroni (1926), de Finetti
(1931), Gini (1938), Dodd (1940) and Aczél (1948). More recent
contributions are the works of Wimp (1986), Hutník (2006), Matkowski
(1999;2002) Jarczyk and Matkowski (2000), J. Marichal (2000), Daróczy and
Hajdu (2005), and Abrahamovic {\em et al.} (2006). Quasi-arithmetic means,
in particular, have been applied in several disciplines. Their functional
form has been used in the theory of copulas under the name of Archimedean
copulas (Genest and MacKay, 1986) and a rich literature can be found under
this name. In the theory of aggregation operators and fuzzy measures, a
growing literature related to the use of quasi-arithmetics includes the
works of Frank (1979), Hajék (1998), Kolesaróvá (2001), Klement {\em et
al.} (1999), Grabish (1995)  and Calvo and Mesiar (2001).

Despite the extensive quasi-arithmetic means literature, to the best of
our knowledge, there is no published work relating quasi-arithmetic means
with covariance functions, whose properties have been extensively studied
both in mathematical analysis and statistical fields. In particular, the
study of covariance functions is intimately related to that of positive
definite functions, the latter being the subject of a considerable
literature in a variety of fields, such as mathematical analysis, abelian
semigroup theory, spatial statistics and geostatistics. For basic facts
about positive definite functions, we refer the interested reader to Berg
and Forst (1975) and to Berg {\em et al.} (1984). The importance of
positive definite functions in the determination of permissible covariance
functions (ordinary and generalised) in spatial statistics was studied in
detail by Christakos (1984). Subsequent considerations, in a spatial and a
spatiotemporal context, include the works of Christakos (1990, 1991,
1992), Sasvári (1994) and Gneiting (1997), among others.

Fundamental properties of covariance functions may be inferred by studying
collections of them considered as convex cones closed in the topology of
point-wise convergence. In this paper, quasi-arithmetic averages and
positive definite functions are combined to gain valuable insight
concerning certain covariance properties. Also, we apply the general
results obtained by our analysis based on the concept of
quasi-arithmeticity to build new classes of stationary and nonstationary
space-time covariance functions. In such a context, we seek answers to the
following questions:
\begin{description} \item[(1)] Consider an arbitrary number $n \in
\mathbb{N}$ of covariance functions, not necessarily defined in the same
space. Their arithmetic average and product ({\em i.e., the geometric
product raised to the $n$-th power}) are valid covariance functions, and
so is the $k$-th power average of covariance functions ($k$ is a natural
positive number). But, what about other types of averages? Since
quasi-arithmetic means constitute a general group that includes the
arithmetic, geometric, power and logarithmic means as special cases, it
seems natural to use quasi-arithmetic representations in order to derive
positive definiteness conditions for quasi-arithmetic averages of
covariance functions. \item[(2)] Can we use the properties of
quasi-arithmetic means to establish important inequalities, ordering and
preference relations and minimal and maximal elements within the class of
covariance functions? \item[(3)] Is it possible to find a class of link
functions that, when applied to the $k$ covariances, can generate valid
nonseparable covariances? If this is the case, other potentially desirable
properties should be examined, e.g., this approach should be as general as
possible and include {\em famous} constructions and separability as
special cases; and it should preserve the margins.
\end{description}

In view of the above considerations, the paper is organised as follows:
Section 2 discusses the background, notation and proposed methodology; it
also provides a very brief introduction to positive definite functions. In
Section 3, the main theoretical results are presented. In particular, we
propose permissibility criteria for the quasi-arithmetic mean of
covariance functions. We derive important covariance inequalities of
Jensen's type as well as ordering and preference relations between
covariance functions. Finally, minimal and maximal elements of the
quasi-arithmetic covariance class are identified, and an associativity
property of this class is provided. In Section 4, we apply our results to
the construction of new families of space-time nonseparable stationary
covariance functions. Also, we extend Stein (2005b) result to a more
general class of spatial nonstationary covariance functions. Several
examples of stationary and nonstationary covariances are proposed and
their mathematical properties discussed. Finally, we study the properties
of quasi-arithmetic random fields (defined later in the paper) in terms of
mean square differentiability and variance. Section 5 concludes with a
critical discussion of the preceding analysis. All the proofs of the
original results derived in this paper can be found in the Appendix.
\section{Background and Methodology}
\subsection{Covariance functions: characterisation and basic properties}

In this section, we shall assume the real mapping ${\rm C}: \mathbb{X}
\times \mathbb{X} \subseteq \RR^d \times \RR^d \to \RR$ to be continuous
and Lebesgue measurable on the domain $\mathbb{X} \times \mathbb{X}$,
where $\mathbb{X}$ can be either a compact space or the entire
$d$-dimensional Euclidean space, $d \in \mathbb{N}$. ${\rm C}$ is a
covariance function if and only if it is positive definite, that is
\begin{equation}
\sum_{i=1}^n \sum_{i=1}^n c_i c_j {\rm C}(\xx_{i},\xx_{j}) \geq 0
\label{permissibility}
\end{equation} for any finite set of real coefficients $\{c_i\}_{i=1}^n$, and for $\xx_1, \ldots,
\xx_n \in \mathbb{X}$. Christakos (1984) calls the covariance condition
(\ref{permissibility}) {\em permissibility}, and throughout the paper we
shall use both this term and that of positive definiteness to characterise
valid covariance functions. A subclass of positive definite functions,
called {\em stationary}, is obtained if $${\rm
C}(\xx_i,\xx_j):=\tilde{{\rm C}}_0(\xx),
$$ with $\xx=\xx_i-\xx_j$ and $\tilde{{\rm C}}_0: \mathbb{X} \to{\RR}$ positive definite
and such that $\tilde{{\rm C}}_0(\boldsymbol{0}) < \infty$.

As shown by Bochner (1933), condition (\ref{permissibility}) is then
equivalent to the requirement that $\tilde{{\rm C}}_0$ is the Fourier
transform $\mathcal{F}$ of a positive bounded measure $\hat{\rm C}_0$ with
support in $\RR^d$, that is
$$\tilde{{\rm C}}_0(\xx):=\mathcal{F}[\hat{\rm C}_0](\xx)= \int_{\RR^d} \ee^{i \boldsymbol{\omega}^{\prime} \xx} \d \hat{\rm C}_0(\boldsymbol{\omega}).
$$ Additionally, if $\hat{\rm C}_0$ is absolutely continuous with respect to the
Lebesgue measure (ensured if $\tilde{\rm C}_0 \in L_1(\mathbb{X})$), then
the expression above can be written as a function of $\d \hat{\rm
C}_0(\boldsymbol{\omega}) = \hat{\rm c}_0(\boldsymbol{\omega}) \d
\boldsymbol{\omega}$, where $\hat{\rm c}_0$ is called the spectral density
of $\tilde{{\rm C}}_0$. For a detailed mathematical discussion of the
Fourier representation in a spatial-temporal statistics context, we refer
the interested reader to volumes by Yaglom (1987) and Christakos (1992).

In what follows, we shall consider some interesting restrictions on the
general class of stationary covariance functions. But first, we need to
introduce some standard notation for arbitrary partitions and operations
between vectors. In order to manipulate arbitrary decompositions of
nonnegative integer numbers, let us consider the set $\expz = \emptyset
\cup \mathbb{Z}_{+} \cup \mathbb{Z}_{+}^2 \cup \mathbb{Z}_{+}^3 \cup
\ldots$ (disjoint union). An element $\mathbf{d}$ of $\expz$ can be
expressed either as $\mathbf{d} = \emptyset$ or as $\mathbf{d} = (d_1,
d_2, \ldots, d_n)$ if $\mathbf{d} \in \mathbb{Z}_{+}^n$ with $n \geq 1$.
In the latter case we denote by $n(\mathbf{d}) = n$ the dimension of
$\mathbf{d}$ and $|\mathbf{d}| = \sum_{i=1}^n d_i$ the length of
$\mathbf{d}$. Both values are taken to be $0$ whenever $\mathbf{d} =
\emptyset$. For $\mathbf{d}, \mathbf{d}' \in \expz$ we say that
$\mathbf{d} \leq \mathbf{d}'$ if and only if $n(\mathbf{d}) =
n(\mathbf{d}')$ and $d_i \leq d'_i$ for all $i=1, 2, \ldots,
n(\mathbf{d})$. Usual vector operations are possible only between elements
of the same dimension. Vectors with all components equal are denoted in
bold symbols, such as $\mathbf{0}, \mathbf{1}$.

Now, the restrictions considered in this paper are of the following type:
\begin{enumerate}
\item The {\bf isotropic} case:
\begin{equation}
\tilde{{\rm C}}_0(\xx):=\tilde{{\rm C}}_1(\| \xx \|), \quad \xx \in
\mathbb{X}, \label{isotropy}
\end{equation}
that is, {\rm C} is said to be represented by the positive definite
function $\tilde{{\rm C}}_1: \RR_{+} \to \RR$, that is
rotation-translation invariant (or radially symmetric), and where $\|.\|$
denotes the Euclidean norm. This is the most popular case in spatial and
spatiotemporal statistics (Matérn, 1960; Yaglom, 1987, Christakos and
Papanicolau, 2000). \item The {\bf component-wise isotropic} case: Let us
consider the $d$-dimensional space $\RR^d$, and let $\mathbf{d}$ be an
element of $\expz$ such that $|\mathbf{d}| = d$ and $\mathbf{1} \leq
\mathbf{d}$. Thus, one can create opportune partitions of the spatial lag
vector $\xx \in \RR^{d}$ in the following way. If $\mathbf{d} = (d_1, d_2,
\ldots, d_n)$ and $\xx \in \RR^d$ we can always write
$$
\xx = (\xx_1, \xx_2, \ldots, \xx_n) \in \RR^{d_1} \times \RR^{d_2} \times
\cdots \times \RR^{d_n}$$ so that:

($\bf{i}$)- $\tilde{{\rm C}}_0(\xx) = \tilde{{\rm C}}_0( \mathbf{k} )$ for
any $\xx, \mathbf{k} \in \RR^d$ if and only if $\| \xx_i \| = \|
\mathbf{k}_i \|$ for all $i=1,2,\ldots, n$.

($\bf{ii}$)- The resulting covariance admits the representation
\begin{eqnarray}
\tilde{{\rm C}}_0(\xx)&:=&\tilde{{\rm C}}_1(\| \xx_{1} \|,\ldots, \|
\xx_{n} \|)
\nonumber \\
&=& \int_{0}^{\infty} \ldots \int_{0}^{\infty} \prod_{i=1}^{n}
\Omega_{d_{i}}(\| \xx_{i} \| r_{i}) \d F(r_{1},\ldots,r_{n})
\label{shapiro}
\end{eqnarray}
\noindent with $\Omega_{d}(t) = \Gamma(d/2) \left (\frac{2}{t} \right )
\mathcal{J}_{(d-2)/2}(t)$, $\mathcal{J}_{d}(.)$ denoting the Bessel
function of the first kind of order $d$ (Abrahamowitz and Stegun, 1965),
$F$ a $n$-variate distribution function and $\tilde{{\rm C}}_1: \RR^n \to
\RR$ positive definite. Thus, equation (\ref{isotropy}) is a special case
of (\ref{shapiro}) ($|\dd|=d$ and $\dd:=d$, a scalar) and its
corresponding integral representation can be readily obtained. The special
case $\mathbf{d}=\boldsymbol{1}$ and $|\mathbf{d}|=d$ is particularly
interesting in the subsequent sections of this paper, as the function
$(x_1,\ldots,x_n) \mapsto \tilde{{\rm C}}_1(|x_1|,\ldots,|x_n|)$, with
$\tilde{{\rm C}}_1$ positive definite, does not depend on the Euclidean
norm, but on the Manhattan or city block distance, with important
implications in spatial and spatiotemporal statistics as pointed out by
Christakos (2000) and Banerjee (2004).
\end{enumerate}
It is worth noticing that covariance functions of the type (\ref{shapiro})
have a property called reflection symmetry (Lu and Zimmerman, 2005), or
full symmetry (Christakos and Hristopulos, 1998). This means that ${\rm
C}(\xx_1,\ldots,\xx_i,\ldots,\xx_n)={\rm
C}(\xx_1,\ldots,-\xx_i,\ldots,\xx_n)=\ldots={\rm
C}(-\xx_1,\ldots,-\xx_i,\ldots,-\xx_n)$. Whenever no confusion arises, in
the remainder of the paper we shall drop the under- and super-script
denoting a stationary or stationary and isotropic covariance function,
respectively.

As far as the basic properties of covariance functions are concerned,
assume that ${\rm C}_i: \RR^{d_i} \to \RR_+$ ($i=1,\ldots,n$) are
positive, continuous and integrable stationary covariance functions with
$d_i \in \mathbb{Z}_{+}$ and let $\dd=(d_1,\ldots,d_n)$ such that
$|\dd|=d$. It is well known that some mean operators preserve the
permissibility of the resulting structure. In particular, if we assume,
without loss of generality, that the $\theta_i$ are nonnegative weights
summing up to one, the following are then permissible on $\RR^d$:
\begin{description}
\item[1] The arithmetic average $${\rm C}_{\rm
A}(\xx)=\sum_{i=1}^{n} \theta_{i} {\rm C}_{i} (\xx_i),$$ \item[2] the non
weighted geometric average up to a power $n$, $${\rm C}_{\rm G}(\xx)=
\prod_{i=1}^{n} {\rm C}_{i}(\xx_i), $$ \item[3] the $k$-power average ($k
\in \mathbb{Z}_+$), up to a power $k$,
$${\rm C}_{\mu^k}(\xx)=\sum_{i=1}^{n} \theta_i C_{i}^k(\xx_i).$$
\item[4] Scale and power mixtures of covariance functions
(Christakos, 1992),
$${\rm C}(\xx)=\int_{\Theta} {\rm C}(\xx;\theta) \d \mu
(\theta),$$for $\mu$ a positive measure and $\theta \in \Theta \subseteq
\RR^p$, $p \in \mathbb{N}$.
\end{description}

\subsection{The methodology. Quasi-arithmetic multivariate compositions}
Quasi-arithmetic averages have been extensively treated by Hardy et al.
(1934). Our methodology generalises the concept of quasi-arithmetic
averages and introduces new formalisms and notation, since our aim is to
find a class of compositions of covariance functions that satisfies
desirable properties.

Let $\Phi$ be the class of real-valued functions $\varphi$ defined in some
domain $D(\varphi) \subset\RR$, admitting a proper inverse $\varphi^{-1}$
defined in $D(\varphi^{-1}) \subset \RR$ and such that $\varphi
(\varphi^{-1} (t)) =t$ for all $t \in D(\varphi^{-1})$. Also, let
$\Phi_{\rm c}$ and $\Phi_{\rm cm}$ be the subclasses of $\Phi$ obtained by
restricting $\varphi$ to be, respectively, convex or completely monotone
on the positive real line. Let us call a {\em quasi-arithmetic class of
functionals} the class
\begin{equation}\label{quasi-functional}
\mathfrak{Q}:= \left \{\psi: D(\varphi^{-1}) \times \cdots \times
D(\varphi^{-1}) \to \RR : \psi(\mathbf{u})= \varphi \left (\sum_{i=1}^n
\theta_{i} \varphi^{-1} (u_i) \right ), \quad \varphi \in \Phi \right \},
\end{equation}
where $\theta_i$ are nonnegative weights and
$\mathbf{u}=(u_1,\ldots,u_n)^{\prime}$, for $n \geq 2$ positive integer.
Also, we shall call $\mathfrak{Q}_{{\rm c}}$ and $\mathfrak{Q}_{{\rm cm}}$
the subclasses of $\mathfrak{Q}$ when restricting $\varphi$ to belong,
respectively, to $\Phi_{\rm c}$ and $\Phi_{\rm cm}$.

If $\psi \in \mathfrak{Q}$, then we should write $\varphi_{\psi}$ as the
function such that: for any nonnegative vector $\mathbf{u}$,
$\psi(\mathbf{u})=\varphi_{\psi}\left ( \sum_{i=1}^n
\theta_i\varphi^{-1}_{\psi}(u_i) \right )$. For ease of notation, we
simply write $\varphi$ instead of $\varphi_{\psi}$, whenever no confusion
arises.


Next, we introduce a new class of functionals that will be used
extensively throughout the paper.

\begin{defn} {\bf Quasi-arithmetic compositions}
If $f_i: \RR^{d_i} \to \RR_{+}$ such that $\cup_{i}^{n} f_i(\RR^{d_i})
\subset D(\varphi^{-1})$ for some $\varphi \in \Phi$, the {\em
quasi-arithmetic composition of $f_1, f_2, \ldots, f_n$ with generating
function $\psi \in \mathfrak{Q}$ is defined as the functional
\begin{equation}
\mathcal{Q}_{\psi}(f_1, \ldots, f_n)(\xx) = \psi \left( f_1(\xx_1) ,
\ldots f_n(\xx_n) \right) \label{eq:composition}
\end{equation}
for $\xx=(\xx_1, \ldots, \xx_n)^{\prime}$, $\xx_i \in \RR^{d_i}$
$\dd=(d_1, \ldots, d_n)^{\prime}$ and $|\dd|=d$.}
\end{defn}
Throughout the paper, we refer to $\psi \in \mathfrak{Q}$ or the
corresponding $\varphi \in \Phi$ as the generating functions of
$\mathcal{Q}_{\psi}$. Note that $\mathcal{Q}_{\psi}(f, \dots, f) = f$ for
any function $f$ and generating function $\psi$.

\begin{exam} \label{example1}
In Table \ref{tab01} we present four basic examples of the
quasi-arithmetic composition considered above. Some conventions are needed
in order to deal with possibly ill-defined values.

\begin{table}
\caption{\label{tab01}Examples of quasi-arithmetic compositions for some
possible choices of the generating function $\varphi \in \Phi$.}
\centering \noindent \begin{tabular}{c} $\begin{array}{|c|c|c|l|} \hline
\multicolumn{1}{|c|}{\varphi(t)}  & \multicolumn{1}{|c|}{\varphi^{-1}(t)}
& \multicolumn{1}{|c|}{\mathcal{Q}_{\psi}(f_1, f_2) (\boldsymbol{\omega})}
&
\multicolumn{1}{|c|}{\text{Remarks}} \\
\hline \hline
             &                  &
     & f_i: \RR^{d_i} \to [0, \infty) \\
\exp(-t)    & - {\rm ln} t         & \prod_{i=1}^n f_i(\xx_i)^{\theta_i}
    & {\rm ln} 0 = -\infty \\
             &                  &
     & \exp(-\infty) = 0 \\
            &                  &
     & \sum \theta_i = 1 \\
\hline
             &                  &
     & f_i: \RR^{d_i} \to [0, \infty) \\
1/t         &  1/t             & \frac{\sum_{i=1}^n \theta_i}{\sum_{i=1}^{n} \frac{\theta_i}{f_i(\xx_i)}} & 1/0 = \infty,\ 1/\infty = 0 \\
             &                  &
     & 0/0 = 0 \\
                 &                  &
     & \sum \theta_i=1 \\
\hline
 & &
    & f_i: \RR^{d_i} \to [0, M] \\
M(1- t/M)_{+} & M(1- t/M)_{+}  &  \sum_{i=1}^n \theta_i f_i(\xx_i)& \text{ for some } M>0 \\
             &                  &
     & (u)_{+} = \max(u, 0)\\
                 &                  &
     & \sum \theta_i = 1\\
\hline - {\rm ln} t    &  \exp(-t)        & \displaystyle{- {\rm ln} \left
(
\sum_{i=1}^n \exp(-\theta_i f_i (\xx_i)) \right )} & f_1, f_2: \RR^d \to \RR \\
\hline
\end{array}$
\end{tabular}
\end{table}
\end{exam}

Ordering relations as well as a minimal element can be found among the set
of quasi-arithmetic compositions of $n \in \mathbb{N}$ fixed functions
indexed by convex generating functions (a maximal element can also be
found only when both fixed functions are upper bounded). From now on, we
shall write $$\mathcal{Q}_{\rm G}(f_1,\ldots,f_n)=\prod_{i=1}^n
f_i^{\theta_i},$$ which is the quasi-arithmetic composition associated
with $\varphi(t)=\exp(-t)$ that generates a geometric average (constraints
on the weights are specified in Table \ref{tab01}). Also, let
$$\mathcal{Q}_{\rm A}(f_1, \ldots, f_n)=\sum_{i=1}^{n} \theta_i f_i,$$
for $f_1, \ldots, f_n: \RR^d \to [0, M]$), be the quasi-arithmetic
composition associated with $\varphi(t)=M(1-t/M)_{+}$ that
 generates the arithmetic average.
Finally, we shall call $$\mathcal{Q}_{\rm
H}(f_1,\ldots,f_n)=\frac{\sum_{i=1}^n \theta_i}{\sum_{i=1}^{n}
\frac{\theta_i}{f_i(\xx_i)}},$$ the quasi-arithmetic composition
associated with $\varphi(t)=1/t$ and generating the harmonic average.

We shall write $g_1 \leq g_2$ whenever $g_1(\xx_1,\ldots,\xx_n) \leq
g_2(\xx_1,\ldots,\xx_n)$ for all $\xx=(\xx_1,\ldots,\xx_n)^{\prime} \in
\RR^d$. Finally, recall that a function $g$ is subadditive whenever
$g(a+b) \leq g(a) + g(b)$ for all $a,b$ in its domain. Thus, we have the
following result:

\begin{prop} \label{prop:archim-ordering}
Given any set of functions $f_1, \ldots, f_n$ and generating functions
$\varphi, \varphi_1, \varphi_2 \in \Phi_{{\rm cm}}$ (which $\psi, \psi_1,
\psi_2 \in \mathfrak{Q}_{{\rm cm}}$ are associated with, respectively),
and given the same set of weights $\{\theta_i\}_{i=1}^n$ for any pairwise
comparison, we have the following point-wise order relations,
\begin{itemize}
    \item $\mathcal{Q}_{\rm G}(f_1,\ldots,f_n) \leq \mathcal{Q}_{\psi}(f_1, \ldots,f_n) \leq
    \sum_{i=1}^n \theta_i f_i$;
    \item $\mathcal{Q}_{\psi}(f_1, \ldots,f_n) \leq
    \mathcal{Q}_{\rm A}(f_1,\ldots,f_n)$ whenever $f_1, \ldots,f_n: \RR^{d_i} \to [0, M]$;
    \item $\mathcal{Q}_{\psi_1}(f_1, \ldots,f_n)  \leq \mathcal{Q}_{\psi_2}(f_1, \ldots,f_n)$ if and only if the function $\varphi_1^{-1} \circ \varphi_2$ is
    subadditive.
\end{itemize}
\end{prop}
This result extends to the functional case those reported by Nelsen
(1999). Its proof can be found in the Appendix.

\subsection{Other useful notions and notation}

A real mapping $\gamma: \mathbb{X} \subseteq \RR^d \to \RR$ is called an
intrinsically stationary variogram (Matheron, 1965) if it is conditionally
negative definite, that is
$$ \sum_{i=1}^{n} \sum_{i=j}^{n}  a_i a_j \gamma(\xx_i-\xx_j) \leq 0 $$ for all
finite collections of real weights $a_i$ summing up to zero and all points
$\xx_i \in \mathbb{X}$. The restriction to the isotropic case is
analogical to that of covariance functions, that is
$\gamma(\xx):=\tilde{\gamma}(\|\xx\|)$, $\xx \in \mathbb{X}$ and
$\tilde{\gamma}: \RR_+ \to \RR$ conditionally negative definite.

A completely monotone function $\varphi$ is a positive function defined on
the positive real line and satisfying $$(-1)^n \varphi^{(n)}(t) \geq 0,
\qquad t>0,
$$ for all $n \in \mathbb{N}$. Completely monotone functions are
characterised in Bernstein's theorem (see Feller, 1966, p.439) as the
Laplace transforms of positive and bounded measures. By a theorem of
Schoenberg (1938), a function ${\rm C}$ is radially symmetric and positive
definite on any $d$-dimensional Euclidean space $\RR^d$ if and only if
${\rm C}(\xx):= \varphi(\|\xx\|^2)$, $\xx \in \RR^d$, with $\varphi$
completely monotonic on the positive real line.

Bernstein functions are positive functions defined on the positive real
line, whose first derivative is completely monotonic. Once again, an
intimate connection with (negative) definiteness arises, as $\gamma$ is a
radially symmetric and conditionally definite function  on any
$d$-dimensional Euclidean space $\RR^d$ if and only if
$\gamma(\xx):=\mathcal{B}(\|\xx\|^2)$, with $\mathcal{B}$ a Bernstein
function.

Sufficient conditions for positive definiteness are stated in Pólya's
criteria (Berg and Forst, 1975) in $\RR^1$; and in Christakos (1984, 1992)
and Gneiting (2001), who extend criteria of the Pólya type to $\RR^{d}$.

\section{Theoretical Results}

In this section we present theoretical results in a general setting, {\em
e.g.} working with arbitrary partitions of $d$-dimensional spaces as
explained previously. In particular, we shall obtain permissibility
criteria for quasi-arithmetic averages of covariance functions on $\RR^d$.
This will be done for ({\bf a}) a general case in which the respective
arguments of the covariance functions used for the quasi-arithmetic
average have no restrictions; ({\bf b}) the restriction to the
component-wise isotropic case; and ({\bf c}) the further restriction to
covariances that are isotropic and defined on the real line. Also, we
shall show some properties of this construction. In particular, we refer
to the associativity of quasi-arithmetic functionals and to the extension
of ordering relations in Proposition \ref{prop:archim-ordering} to the
case of compositions of covariance functions. The proof of these new
results can be found in the Appendix.
\begin{prop} \label{christakos1}

\begin{description}
\item[(a) General case] Let $\varphi \in \Phi_{{\rm cm}}$ and
${\rm C}_{i}: \RR^{d_i} \to \RR$ ($i=1,\ldots,n$) be continuous stationary
covariance functions such that $\cup_{i}^{n} {\rm C}_{i}(\RR^{d_i})
\subset D(\varphi^{-1})$ and ${\bf d}=(d_1,\ldots,d_n)^{\prime}$,
$|\dd|=d$.  If the functions $\xx_i \mapsto \varphi^{-1} \circ
C_{i}(\xx_i)$, $i=1,\ldots,n$, are intrinsically stationary variograms on
$\RR^{d_i}$, then
\begin{equation} \label{eq:christakos1}
\mathcal{Q}_{\psi} \left ( {\rm C}_1, \ldots, {\rm C}_n \right )
(\xx_1,\ldots,\xx_n)
\end{equation} is a stationary covariance function on $\RR^d$.
\item[(b) Component-wise isotropy] Let $\varphi \in \Phi_{{\rm
cm}}$ and ${\rm C}_{i}: \RR^{d_i} \to \RR$ ($i=1,\ldots,n$) be continuous
stationary and {\bf isotropic} covariance functions such that
$\cup_{i}^{n} {\rm C}_{i}(\RR^{d_i}) \subset D(\varphi^{-1})$ and ${\bf
d}=(d_1,\ldots,d_n)^{\prime}$, $|\dd|=d$. If the functions $x \mapsto
\varphi^{-1}\circ {\rm C}_i(x)$ are Bernstein functions on the positive
real line, then
\begin{equation}\label{eq:christakos2}  \mathcal{Q}_{\psi}\left (
{\rm C}_1, \ldots, {\rm C}_n \right ) (\| \xx_1 \|,\ldots,\| \xx_n \|)
\end{equation} is a stationary and fully symmetric covariance function on $\RR^d$.
\item[(c) Univariate covariances] Let $\varphi \in \Phi_{{\rm
cm}}$ and ${\rm C}_{i}: \RR \to \RR$ ($i=1,\ldots,n$) be continuous
stationary and {\bf isotropic} covariance functions defined on the real
line such that $\cup_{i}^{n} {\rm C}_{i}(\RR) \subset D(\varphi^{-1})$ and
$|\dd|=n$. If the functions $x \mapsto \varphi^{-1}\circ {\rm C}_i(x)$ are
continuous, increasing and concave on the positive real line, then
\begin{equation} \label{eq:christakos3} \mathcal{Q}_{\psi} \left (
{\rm C}_1, \ldots, {\rm C}_n \right ) (|x_1|,\ldots,|x_n|)
\end{equation} is a stationary and fully symmetric covariance function on $\RR^n$.
\end{description}
\end{prop}

It is worth noticing that case (\ref{eq:christakos3}) represents a
covariance permissibility condition that does not depend on the Euclidean
metric, as it is function of the Manhattan or city block distance. For a
detailed discussion of the limitations of the Euclidean norm-dependent
covariances, see Banerjee (2004).

The previous result is of particular importance, since it has implications
both in terms of averages of covariance functions and mixtures thereof. In
the following we derive sufficient permissibility conditions for geometric
and harmonic averages of covariance functions.

\begin{cor} \label{coro1} ({\bf Geometric average})
Let $C_i: \RR^{d_i} \to \RR$ ($i=1,\ldots,n$) be continuous permissible
covariance functions. Let $\theta_i$ ($i=1 \ldots, n$) be nonnegative
weights summing up to one. If the functions $x \mapsto -{\rm ln} \left
(C_i(x) \right )$, $x>0$ satisfy any of the relevant conditions described
in {\bf (a)},{\bf (b)} and {\bf (c)} of Proposition \ref{christakos1}
above, then
$$\mathcal{Q}_{\rm G}(C_1,\ldots,C_n)=\prod_{i=1}^{n}
C_{i}^{\theta_i} $$ is a covariance function.
\end{cor}
An example of this setting can be found by using the function $x \mapsto
(1+x^{\delta})^{-\varepsilon}$, $x$ positive argument, $\delta \in (0,2]$
and $\varepsilon$ positive, also known as generalised Cauchy class
(Gneiting and Schlather, 2004). One can verify that the composition of
this function with the natural logarithm is continuous, increasing and
concave on the positive real line for $\delta \in (0,1]$. Another function
satisfying these requirements is the function $x \mapsto
\exp(-x^{\delta})$, which is completely monotonic for $\delta \in (0,1]$.
It should be stressed that these permissibility criteria do not apply to
compactly supported covariance functions, such as the spherical model
(Christakos, 1992).

\begin{cor} ({\bf Harmonic average}) Let $C_i: \RR^{d_i} \to
\RR$ ($i=1,\ldots,n$) be continuous permissible covariance functions. If
the functions $x \mapsto C_i(x)^{-1}$, $x>0$ satisfy any of the relevant
conditions described in {\bf (a)},{\bf(b)} and {\bf (c)} of Proposition
\ref{christakos1} above, then for $\theta_i \geq 0$ such that $\sum_i
\theta_i=1$, the
$$\mathcal{Q}_{\rm H}(C_1,\ldots,C_n)=\frac{\sum
\theta_i}{\sum \frac{\theta_i}{C_i}}
$$ is a covariance function.
\end{cor}

In order to complete the picture about quasi-arithmetic covariance
functions, it would be desirable to establish at least some of their
algebraic properties.  The following results show some important features
of the theoretical construction obtained from Proposition
\ref{christakos1}.

\begin{prop} \label{christakos4} ({\bf Associativity}).
Consider the same arrangements as in Proposition \ref{christakos1} and set
$\theta_i=1/n$ ($i=1,\ldots,n$). Let $\varphi_1 \in \Phi_{{\rm cm}}$ and
$\varphi_2 \in \Phi$. If the functions $x \mapsto \varphi^{-1} \circ
\varphi_2(x)$, $x \mapsto \varphi^{-1}_{2} \circ {\rm C}_{i}(x)$ (
$i=1,\ldots,k$) and $x \mapsto \varphi_1^{-1} \circ C_j(x)$
($j=k+1,\ldots,n$; $t>0$ and $k<n$) satisfy any of the relevant conditions
described in {\bf (a)},{\bf(b)} and {\bf (c)} of Proposition
\ref{christakos1} above, then
\begin{equation} \label{eq:associativity0}
\mathcal{Q}_{\psi_{1}} \left ( \mathcal{Q}_{\psi_2}({\rm
C}_{1},\ldots,{\rm C}_{k}) , {\rm C}_{k+1},\ldots, {\rm C}_n \right
)(\xx_1,\ldots,\xx_k,\xx_{k+1},\ldots,\xx_n),
\end{equation}
\begin{equation} \label{eq:associativity}
\mathcal{Q}_{\psi_{1}} \left ( \mathcal{Q}_{\psi_2}({\rm
C}_{1},\ldots,{\rm C}_{k}) , {\rm C}_{k+1},\ldots, {\rm C}_n \right
)(\|\xx_1\|,\ldots,\|\xx_k\|,\|\xx_{k+1}\|,\ldots,\|\xx_n\|)
\end{equation}  and
\begin{equation} \label{eq:associativity2}
\mathcal{Q}_{\psi_{1}} \left ( \mathcal{Q}_{\psi_2}({\rm
C}_{1},\ldots,{\rm C}_{k}) , {\rm C}_{k+1},\ldots, {\rm C}_n \right
)(|x_1|,\ldots,|x_k|,|x_{k+1}|,\ldots,|x_n|)
\end{equation} and accordingly any coherent permutation of $\varphi_1,\varphi_2$ with
${\rm C}_{i}$ ($i=1,\ldots,n$) are covariance functions.
\end{prop}
One may notice that we deviate from the associativity condition defined in
Bemporad (1926) and (in a more general form called decomposability) in
Marichal (2000). Nevertheless, it can be readily verified that the
associativity condition of Bemporad is always satisfied for our
construction, whereas the strong and weak decomposability of Marichal
(2000) is certainly satisfied by construction (\ref{eq:christakos3}), but
not, in general, by constructions (\ref{eq:christakos1}) and
(\ref{eq:christakos2}).

\begin{prop} \label{prop:archim-ordering2} ({\bf Ordering and preference
relations}) For any set of covariance functions ${\rm C}_1, \ldots, {\rm
C}_n$ and any arbitrary generating functions $\varphi, \varphi_1,
\varphi_2 \in \Phi_{{\rm cm}}$ with associated $\psi, \psi_1, \psi_2 \in
\mathfrak{Q}_{{\rm cm}}$, if ($i$) the constraints in either {\bf (a)},
{\bf (b)} or {\bf (c)} of Proposition \ref{christakos1} are satisfied and
($ii$) the same set of weights $\{\theta_i\}_{i=1}^k$ is used for any
pairwise comparison, then we have the point-wise order relations,
\begin{itemize}
    \item[(1)] $\mathcal{Q}_{\rm G}({\rm C}_1,\ldots,{\rm C}_n) \leq \mathcal{Q}_{\psi}({\rm C}_1, \ldots,{\rm C}_n) \leq
    \sum_{i=1}^n \theta_i {\rm C}_i$;
    \item[(2)] $\mathcal{Q}_{\rm H}({\rm C}_1, \ldots,{\rm C}_n) \leq
    \mathcal{Q}_{\rm A}({\rm C}_1,\ldots,{\rm C}_n)$ whenever ${\rm C}_1, \ldots,{\rm C}_n: \RR^{d_i} \to [0, M]$;
    \item[(3)] $\mathcal{Q}_{\psi_1}({\rm C}_1, \ldots,{\rm C}_n)  \leq \mathcal{Q}_{\psi_2}({\rm C}_1, \ldots,{\rm C}_n)$ if and only if the function $\varphi_1^{-1} \circ \varphi_2$ is
    subadditive.
\end{itemize}
\end{prop}
The proof is omitted as it follows the same arguments as that of
Proposition \ref{prop:archim-ordering}. It should be stressed that the
first inequality shows that quasi-arithmetic covariances satisfy Jensen's
inequality. This fact has some implications on the variance of the
stationary random field generated by a quasi-arithmetic operator. This
will be discussed in the subsequent sections. Also, it is important to
specify that other results, in particular inequalities involving
quasi-arithmetic operators, can be readily extended to the case of
quasi-arithmetic compositions of covariance functions. This is the case,
for instance, of Theorem 1 in Abramovich {\em et al.} (2004).

\section{Using quasi-arithmetic functionals in the construction of nonseparable space-time covariance functions}
\subsection{Review of space-time covariance functions}

Natural (physical, health, cultural etc.) systems involve various
attributes, such as atmospheric pollutant concentrations, precipitation
fields, income distributions, and mortality fields. These attributes  are
characterised by spatial-temporal variability and uncertainty that may be
due to epistemic and ontologic factors. In view of the prohibiting costs
of spatially dense monitoring networks, one often aims to develop a
mathematical model of the natural system in a continuous space-time
domain, based on sequential observations at a limited number of monitoring
stations. This kind of problem has been a motivation for the development
of the {\em spatio-temporal random field} ($S/TRF$) theory; see Christakos
(1990, 1991, 1992) and Christakos and Hristopulos (1998) for a detailed
discussion of the ordinary and generalised $S/TRF$ theory and its various
applications. In the following, we slightly change notation in order to be
consistent with classical nomenclature in the Geostatistical literature
(Cressie and Huang, 1999; Gneiting, 2002). Let $\{Z(\ss,t), (\ss,t) \in
\RR^d \times \RR \}$ be a real-valued $S/TRF$, where $\ss,t$ denote
respectively the spatial and temporal position. Then, the function
 $\textrm{C}_{s,t}(\ss_1,t_1,\ss_2,t_2)=\textrm{cov}(Z(\ss_{1},t_{1}),Z(\ss_{2},t_{2}))$,
 defined on the product space $\RR^d \times \RR \times \RR^d \times \RR$,
 is the covariance function of the associated $S/TRF$ $Z$ if and only if it is
 permissible, {\em i.e.} satisfies (\ref{permissibility}).
When referring to the spatial index, the term {\em homogeneity} instead of
weak stationarity is equivalently adopted (Christakos, 1990, 1992). Thus,
under the assumption of spatial homogeneity and temporal stationarity
(sometimes, simply called spatio-temporal stationarity in the weak sense),
the underlying $S/TRF$, ergo denoted by $H/S$, has finite and constant
mean and the covariance function, defined on the product space $\RR^d
\times \RR$, is such that
$\textrm{cov}(Z(\ss_{1},t_{1}),Z(\ss_{2},t_{2}))={\rm C}_{s,t}(\hh,u)$,
with $(\hh,u) = (\ss_{i} - \ss_{j},t_{i} - t_{j}) \in \RR^{d} \times \RR$
denoting the spatio-temporal separation
 vector, and the $\hh$ and $u$ denoting the spatial and the temporal lags, respectively.
The special case of isotropy in the spatial component and symmetry in the
the temporal one is denoted as
\begin{equation}
\textrm{C}_{s,t}(\hh,u):=\widetilde{\textrm{C}}_{s,t}(\| \hh \|, |u|),
\label{full_sym}
\end{equation}
\noindent where the $\|.\|$ denotes the Euclidean norm. Obviously,
(\ref{full_sym}) is fully symmetric.

Another popular assumptions concerning the $S/TRF$ model is that of
separability, that is (Christakos and Hristopulos, 1998; Mitchell \emph{et
al.}, 2004) \begin{equation}
\textrm{C}_{s,t}(\hh,u)=\textrm{C}_{s,t}(\hh,0)\textrm{C}_{s,t}(\textbf{0},u).
\label{separable}
\end{equation}
\noindent In other words, separability means that the spatio-temporal
covariance structure factors into a purely spatial and a purely temporal
component, which allows for computationally efficient estimation and
inference. Consequently, separable covariance models have been popular
even in situations in which they are not physically justifiable. Another
interesting aspect is that separable covariances are also fully symmetric,
whereas the converse is not necessarily true.

It has been argued in the relevant literature that separable models allow
for ease of computation and dimensionality reduction, as the space-time
covariance matrix is obtained through the Kronecker product of the
marginal spatial and temporal ones. However, separability is an
unrealistic assumption for many applications, since it implies a
considerable loss of information about important interactions between the
spatial and temporal variations. Therefore, various techniques have been
introduced for generating different classes of nonseparable
spatio-temporal covariance models.  Most of these techniques have been
developed in the context of applied {\em stochastics} analysis and include
H/S as well as non-H/S covariances (e.g., Yaglom, 1948; Gandin and
Boltenkov, 1967; Monin and Yaglom, 1965, 67; Christakos, 1990, 1991, 1992,
2000; Christakos and Hristopulos, 1998; Christakos {\em et al.}, 2002,
2005; Kolovos {\em et al.}, 2004).  Also, covariance models have been
developed in the context of spatio-temporal {\em statistics} (e.g., Jones
and Zhang, 1997; Cressie and Huang, 1999; De Cesare \textit{et al.}, 2001;
Fuentes and Smith, 2001; Gneiting, 2002; Ma, 2003; Fernández-Casal {\em et
al.}, 2003; Stein, 2005a; Porcu {\em et al.}, 2006a).

Being the Laplace transform of positive bounded measures, completely
monotone functions are particularly appealing for the construction of
stationary and nonstationary space-time covariances. In particular, they
are intimately connected with the concept of {\em mixture-based}
covariance functions, that has been repeatedly used by several authors. In
the stationary case, see Gneiting (2002), Ma (2002;2003), Ehm {\em et al.}
(2003), Fernández-Casal {\em et al.} (2003), Porcu {\em et al.} (2006a).
In the nonstationary case, Paciorek and Schervish (2004;2006), Stein
(2005b), and Porcu {\em et al.} (2006b) have made use of this technique.

Also, mixture-based covariances have been developed with less
sophisticated instruments than completely monotone functions. This is the
case of the so-called product sum model (De Cesare, 2001) and their
extensions (De Iaco {\em et al.}, 2002; Porcu {\em et al.}, 2006c). This
group of authors builds nonseparable space-time covariances through simple
application of the basic properties of covariances seen as a convex cone.
The mixture-based procedures and the basic properties of covariance
functions are properly combined in the present paper.

In the following, a stationary RF with a quasi-arithmetic covariance
function will be called {\em quasi-arithmetic random field} and denoted
with the acronym QARF.
\subsection{On the representation and smoothness properties of QARF}

In this section we focus on the representation of a QARF and discuss the
smoothness properties in terms of (mean square) partial differentiability.
These properties are intimately related to those of the associated
covariance function.

Let  $Z_{i}(s)$ be univariate mutually independent continuous weakly
stationary Gaussian random processes defined on the real line ($i=1,
\ldots, (d+1)$; $s \in \RR$; and $d \in \mathbb{Z}_{+}$). In particular,
let the process $Z_{d+1}$ be continuously indexed by time $t$. Consider
also a ($d+1$)-dimensional nonnegative random vector
$\boldsymbol{R}=(R_1,\ldots,R_{d+1})^{'}$ with $R_i$ independent of
$Z_{i}$. Let the univariate covariances ${\rm C}_{s_{i}}$ and the temporal
covariance ${\rm C}_t$ be respectively associated with $Z_i$,
$i=1,\ldots,d$ and $Z_{d+1}$. In the following we shall assume these
covariances to be stationary, symmetric, and of the type ${\rm
C}_{s_{i}}(h_i)=\exp (-\nu_i(|h_i|))$, $i=1,\ldots,d$, and $C_t(u)=\exp
(-\nu_t(|u|))$, with $\nu_i=\varphi^{-1} \circ {\rm C_i},
\nu_t=\varphi^{-1} \circ {\rm C_t}$, where the $\varphi \in \Phi_{\rm cm}$
and ${\rm C}_i$ are positive definite and such that the compositions
$\nu_i$ are continuous, increasing and concave on the positive real line.
Positive definiteness of this construction is guaranteed by direct
application of the theorem of Schoenberg (1938) and according to a P{\'
o}lya type criterion (see Berg and Forst (1975), Proposition 10.6).

We are interested in inspecting the properties of the following stationary
spatio-temporal scale mixture-based random field, defined on $\RR^{d}
\times \RR$, \begin{equation} Z(\ss,t)= Z_{d+1}(t R_{d+1}) \prod_{i=1}^{d}
Z_{i}(s_iR_i), \label{representation}
\end{equation} with $\ss=(s_1,\ldots,s_d)^{'} \in \RR^d$ and $t \in \RR$. It can
be easily seen that the covariance structure associated to this random
field is nonseparable, as
\begin{equation}
{\rm C}_{s,t}(\hh,u)= \int_{\RR_{+}^{d+1}} \exp \left
(-\sum_{i=1}^{d}\nu_{i}(|h_i|)r_{i}-\nu_t(|u|)r_{d+1} \right ) \d F({
\rr}), \label{representation2}
\end{equation}
with $\hh=(h_1, \ldots, h_d)^{'} \in \RR^d$, $u \in \RR$ and $\rr=(r_1,
\ldots, r_{d+1})^{'} \in \RR^{d+1}$, and where $F$ is the ($d+1$)-variate
distribution function associated to the random vector $\boldsymbol{R}$. If
$F$ is absolutely continuous with respect to the Lebesgue measure, then
previous representation can be reformulated with respect to the
($d+1$)-variate density, say $f$, that is
$${\rm C}_{s,t}(\hh,u)= \int_{\RR_{+}^{d+1}}
\exp \left (-\sum_{i=1}^{d}\nu_{i}(|h_i|)r_{i}-\nu_t(|u|)r_{d+1} \right )
f(\rr) \d \rr.
$$

It can be seen easily that this construction allows for the case of
separability if and only if the integrating ($d+1$)-dimensional measure
$F$ (or equivalently its associated density $f$) factorises into the
product of ($d+1$) marginal ones, {\em i.e.} if the nonnegative random
vector $\boldsymbol{R}$ has mutually independent components.

Now, if we suppose that the measure $F$ is concentrated on the line
$r_1=\ldots=r_{d+1}=r$ and that $\varphi^{-1}$ is such that
$\varphi:=\mathcal{L}[F]$, {\em i.e.} the Laplace transform of the
positive univariate measure $F$, then we obtain that the QARF is a special
case of (\ref{representation2}). For this random field, further inferences
may be made about its mean square partial differentiability if,
additionally, we assume the function $t \mapsto \exp(-\nu_i(t))$ to be
absolutely integrable on the positive real line ($i=1,\ldots,(d+1)$). In
this case, one can show (technicalities can be found in the appendix) that
the $k$-th order mean square partial derivative with respect to the $i$-th
coordinate exists and is finite whenever the function $\chi(r)=
\int_{[0,\infty)} \omega_i^{2k} \hat{\rm c}_i(\omega_i; r) \d \omega_i$,
with $\hat{\rm c}_i:=\mathcal{F}^{-1}[\exp(- r \nu_i)]=\int_{\RR}\exp(-i
\omega_i h_i-r\nu_i(h_i)) \d h_i$ is measurable with respect to $F$
($i=1,\ldots,d$).

It should be stressed that Proposition \ref{prop:archim-ordering2} gives
some more information about the characteristics of the underlying QARF in
terms of variance, as QARF can be ordered with respect to their minimum or
maximum variance. Thus, the QARF generated by $\psi:=A$ has the largest
variance among all the other QARF, for any choice of $\psi$.

Finally, observe that the trivial quasi-arithmetic composition, obtained
by setting ${\rm C}_i:= \varphi \in \Phi_{\rm cm}$ ($i=1,\ldots,d$)
preserves permissibility and results in a model of the type
$${\rm C}_{s}(\hh)=\varphi \left ( \boldsymbol{\theta}^{\prime} \hh \right ) $$
where $\boldsymbol{\theta}=(\theta_1,\ldots,\theta_d)^{\prime} \in \RR^d$
and $\hh=(\hh_1,\ldots,\hh_d)^{\prime}$ is the arbitrary partition of the
spatial lag vector. Thus, the trivial case results in a composition of a
completely monotonic function with an affine function, and can be used for
modelling geometric anisotropies.

\subsection{Applications of quasi-arithmeticity in the construction of stationary nonseparable space-time covariances}

The results presented in the previous section can be useful in the
construction of space-time covariance functions. For this purpose, some
considerations are in order. The functional in equation
(\ref{quasi-functional}) should be adapted in the spatio-temporal case. It
does not make sense to consider a weighted average of a spatial covariance
with a temporal one. Also, the use of weights forbids, in this case, the
construction of nonseparable model admitting separability as a special
case. For this reason, we suggest to reduce the class in equation
(\ref{quasi-functional}) to the case of trivial weights, {\em e.g.}
$\theta_i=1$, $\forall i$. Then, it is more appropriate to call the class
(\ref{quasi-functional}) {\em Archimedean}, in analogy with the class
built in Genest and MacKay (1986). It should be mentioned that one can
easily prove that the restriction to trivial weights does not affect the
permissibility of the resulting covariance function, provided that either
one of the constraints imposed in cases {\bf (a)}, {\bf (b)} or {\bf (c)}
of Proposition \ref{christakos1} are fulfilled. By analogy with the
construction of copulas, and following Genest and MacKay (1986), the
definition of {\em generator} for $\varphi \in \Phi_{{\rm cm}}$ is more
meaningful.

Working this way, one can obtain some new families of covariance functions
whose analytical expressions are familiar for those interested in
probabilistic modelling through copulas. It should be noticed that all the
families we propose in this section include separable covariance models as
special cases, depending on the parameter values of the generators.

\subsubsection*{Example 1. The Clayton family}
\noindent Consider the completely monotone function
\begin{equation}
\varphi(x)=(1+x)^{-1/\lambda_{1}}, \qquad x>0, \label{example1_a}
\end{equation}
\noindent where $\lambda_{1}$ is a nonnegative parameter with inverse
$\varphi^{-1}(y)=y^{-\lambda_{1}}-1$. Observe that (\ref{example1_a}) is
the generator of the Clayton family of copulas; refer to Genest and MacKay
(1986) for mathematical details about this class.

Also, consider the covariance functions ${\rm C}_{s}(\hh)=(1+\|\hh
\|)^{-1/\lambda_{2}}$ and ${\rm C}_{t}(u)=(1+|u|)^{-1/\lambda_{3}}$, for
$\lambda_2, \lambda_3$ positive parameters. It is easy to verify that
$\varphi^{-1} \left ( {\rm C}_{i}(y)\right
)=(1+y)^{\lambda_{1}/\lambda_{i}}$ ($i=2,3$) is a positive function whose
first derivative is completely monotone if and only if
$\lambda_{1}<\lambda_{i}$. Under this constraint, and applying case {\bf
(b)} of Proposition \ref{christakos1}, we find that
\begin{equation}
{\rm C}_{s,t}(\hh,u)=\mathcal{Q}_{\psi}({\rm C}_{s},{\rm
C}_{t})(\hh,u)=\sigma^2 \left [ (1+ \|\hh \|)^{\rho_{1}}+(1+
|u|)^{\rho_{2}}-1 \right ]^{-1/\lambda_{1}} \label{example1_b}
\end{equation}
\noindent is a valid nonseparable stationary fully symmetric
spatio-temporal covariance function, with
$\rho_{i}=\lambda_{1}/\lambda_{i}$, $\lambda_{i}>0$ ($i=2,3$) and
$\sigma^2$ a nonnegative parameter denoting the variance of the underlying
$S/T$ process. It is interesting that both margins (the spatial and the
temporal one) belong to the generalised Cauchy class. This is desirable
for those interested in the local and global behaviour of the underlying
random field.

Another covariance function that preserves Cauchy margins can be obtained
through the following {\em iter}. Consider the function $x \mapsto
\varphi(x)=x^{-\alpha}$, $t>0$, that belongs to $\Phi_{{\rm cm}}$ for any
positive $\alpha$, being the Laplace transform of the function $x \mapsto
x^{\alpha-1}/\Gamma(\alpha)$, with $\alpha$ positive parameter and
$\Gamma$ the Euler Gamma function. Also, for the spatial margin consider
${\rm C}_{s}(\hh)=(1+\|\hh\|^{\delta})^{-\varepsilon}$ that belongs to
$\Phi_{{\rm cm}}$ if and only if $\delta \in (0,2]$ and $\varepsilon$ is
strictly positive. Finally, let the temporal margin be of the type ${\rm
C}_{t}(u)=|u|^{-\rho}$, which is not a stationary covariance function but
respects the composition criteria, as $\varphi^{-1} \circ {\rm C}_{t}$ is
continuous, increasing and concave on the positive real line if and only
if $\alpha < \rho$. Then, it is easy to prove that
$\mathcal{Q}_{\psi}({\rm C}_{S},{\rm C}_{t})(\hh,u)$ is a valid space-time
function if and only if $\alpha < \varepsilon$ and $\alpha < \rho$, and
that this covariance function has margins of the Cauchy type.

\subsubsection*{Example 2. The Gumbel-Hougard family} Consider the completely monotone function, the
so-called positive stable Laplace transform
\begin{equation}
\varphi(x)=e^{-x^{1/\lambda_{1}}}, \qquad x >0, \label{example2_a}
\end{equation}
\noindent where $\lambda_{1} \geq 1$. Equation (18) admits the inverse
$\varphi^{-1}(y)=(-{\rm ln}(y))^{\lambda_{1}}$. This is the generator of
the Gumbel-Hougard family of copulas, whose mathematical construction and
details are exhaustively described in Nelsen (1999). Consider two
respectively spatial and temporal covariance functions admitting the same
analytical form, {\em i.e.} ${\rm C}_s(\hh)= \exp(-\|\hh\|^{1/\lambda_2})$
and ${\rm C}_t(u)= \exp(-|u|^{1/\lambda_3})$. Then, it can be easily
verified that $\varphi^{-1} \left ( {\rm C}_{s}(y)\right
)=y^{\lambda_{1}/\lambda_{2}}$, and $\varphi^{-1} \left ( {\rm
C}_{t}(y)\right )=y^{\lambda_{1}/\lambda_{3}}$ are always positive for
$y>0$ and possess completely monotone first derivatives if and only if
$\lambda_{1} < \lambda_{i}$, $i=2,3$. So we get that
\begin{equation}
{\rm C}_{s,t}(\hh,u)=\mathcal{Q}_{\psi}({\rm C}_{s},{\rm C}_t)(\hh,u)=
\sigma^2 \exp \left (-(\|\hh \|^{\rho_{1}}+|u|^{\rho_{2}})^{1/\lambda_{1}}
\right ),
\end{equation}
\noindent is a permissible nonseparable stationary fully symmetric
spatio-temporal covariance function, with
$\rho_{i}=\lambda_{1}/\lambda_{i}$, $i=2,3$ and $\sigma^2$ as before.

\subsubsection*{Example 3. The power series family}

\noindent The so-called power series
\begin{equation}
\varphi(x)=1-(1-\exp(-x))^{1/\lambda_{1}}, \qquad x>0, \label{example3_a}
\end{equation}
\noindent with $\lambda_{1} \geq 1$, admits the inverse
$\varphi^{-1}(y)=-{\rm ln}(1-(1-y)^{\lambda_{1}})$. Suppose a spatial and
a temporal covariance function of the same analytical form as in Example
2. The composition $\varphi^{-1} ({\rm C}_{i}(y))=-{\rm
ln}(1-(1-\exp(-y))^{\lambda_{1}/\lambda_{i}})$, $i=2,3$, is always
positive for $y>0$ and admits a completely monotone first derivative if
and only if $\lambda_{1} < \lambda_{i}$. So we get that
\begin{eqnarray}
\mathcal{Q}_{\psi}({\rm C}_{S},{\rm C}_{T})(\hh,u) &=&  1- (1-\exp(-\|\hh
\|))^{\rho_{1}}-(1-\exp(-|u|))^{\rho_{2}}  \nonumber \\
&& +(1-\exp(-\|\hh \|))^{\rho_{1}}(1-\exp(-|u|))^{\rho_{2}}
\end{eqnarray}
\noindent is a permissible nonseparable stationary fully symmetric
spatio-temporal covariance function, with
$\rho_{i}=\lambda_{1}/\lambda_{i}$ ($i=2,3$).

\subsubsection*{Example 4. The semiparametric Frank family}
Here we show that a nonstationary covariance function can be obtained,
starting from the proposed setting, even if the arguments of the
quasi-arithmetic functional are not covariance functions.

 The Frank family of copulas (Genest, 1987) is
generated by the function $\varphi(x)=\frac{1}{\lambda} {\rm ln}\left (
1-(1-\ee^{-\lambda})\ee^{-x}\right )$ with inverse $\varphi^{-1}(y)=-{\rm
ln}\left ( (1-\ee^{-\lambda y}) / (1- \ee^{-\lambda}) \right)$. Nelsen
(1999) shows that, for $\lambda$ positive, $\varphi$ is the composition of
an absolutely monotonic function with a completely monotonic one, {\em
i.e.} a completely monotonic function. As far as the inverse is concerned,
it is easy to see that $\varphi^{-1} \circ \gamma$ ($\gamma$ an
intrinsically stationary variogram) is negative definite. Thus, the $${\rm
C}_{s,t}(\hh,u)= -\frac{1}{\lambda} {\rm} ln \left (1+
\frac{(1-\ee^{-\lambda \gamma_{S}(\hh)})(1-\ee^{-\lambda
\gamma_{T}(u)})}{1-\ee^{-\lambda}} \right ),$$ for $\gamma_{S},
\gamma_{T}$ intrinsically stationary variograms defined on $\RR^d$ and
$\RR$, respectively, is a stationary space-time covariance function.

\subsection{Quasi-arithmeticity and nonstationarity in space
}

The direct construction of spatial covariances that are nonstationary is
anything but a trivial fact. Only few contributions refer to this kind of
construction: see, among them, Christakos (1990, 1991), Christakos and
Hristopoulos (1998) and Kolovos {\em et al.} (2004). More recent
contributions can be found in Stein (2005b) and Paciorek and Schervish
(2006). It seems that something more could be done concerning the
construction of nonstationary covariances, knowing that stationarity is
very often an unrealistic assumption for several physical and natural
processes.

In this section, we show that both approaches discussed in Stein (2005b)
and Paciorek and Schervish (2006), admit a natural extension through the
use of quasi-arithmetic functionals. We need to introduce some more
notation concerning the restriction of the class $\Phi_{\rm cm}$. For
abuse of notation, a completely monotonic function is the Laplace
transform of some positive and bounded measure, so that $\varphi \in
\Phi_{\rm cm}$ if and only if $\varphi := \mathcal{L}[F]$.

%

\begin{thm} \label{nonstationary}
Let $\Sigma$ be a mapping from $\RR^p \times \RR^p$ to positive definite
$p \times p$ matrices, $F$ a nonnegative measure on $\RR_{+}$, $\varphi_1,
\varphi_2 \in \Phi_{\rm cm}$ and $g$ a nonnegative function such that, for
any fixed $\ss \in \RR^p$, $h_ß=\varphi_2^{-1} \circ g(.;\ss) \in
L^{1}(F)$. Define $\Sigma(\ss_1,\ss_2)=1/2(\Sigma(\ss_1)+\Sigma(\ss_2))$
and
$Q(\ss_1,\ss_2)=(\ss_1-\ss_2)^{\prime}\Sigma(\ss_1,\ss_2)^{-1}(\ss_1-\ss_2)$.
Then, \begin{equation} {\rm
C}_{s}(\ss_1,\ss_2)=\frac{|\Sigma(\ss_1)|^{1/4}|\Sigma(\ss_2)|^{1/4}}{|\Sigma(\ss_1,\ss_2)|^{1/2}}\int_{0}^{\infty}
\varphi_1\left (Q(\ss_1,\ss_2)\tau \right )
\mathcal{Q}_{\psi_2}(g_{\ss_1},g_{\ss_2})(\tau) \d F (\tau)
\label{nonstationarity_theorem}
\end{equation}  with $\theta_i=1$, $i=1,2$, is a nonstationary covariance function on $\RR^d \times
\RR^d$.
\end{thm}

Some comments are in order. One can see that Stein's (2005b) result is a
special case of (\ref{nonstationarity_theorem}), under the choice
$\varphi_1(t)=\exp(-t)$, $t$ positive, and $\psi_2:=G$. So is the result
in Paciorek and Schervish (2006). Nevertheless, the form we propose has
some drawbacks that need to be noticed explicitly. The first problem is
that it is very difficult to obtain a closed form for expression
(\ref{nonstationarity_theorem}), unless one chooses $\psi_2:=G$. Secondly,
if one's purpose is to generalise the Matérn covariance function (Matérn,
1960) to the nonstationary case, then the problem has already been solved
by Stein (2005b) and his approach is in our opinion the most suitable, as
he finds a Matérn type covariance that allows for a spatially varying
smoothness parameter and for local geometric anisotropy. On the one hand,
the Matérn covariance possesses some desirable features (i.e., it allows
for arbitrary levels of smoothness of the associated random field).  On
the other hand, there are other covariance functions that are of
considerable interest to the statistical and scientific communities, such
as the nonstationary ones. In particular, we refer to the Cauchy class,
whose properties (in terms of decoupling of the global and local behaviour
of the associated random field) have been discussed in Gneiting and
Schlather (2004). After several trials, we did not succeed in obtaining a
Cauchy type nonstationary covariance. Several examples of covariances can
be derived that belong to the class (\ref{nonstationarity_theorem}) and
can be numerically integrated. Here, we propose some closed form that is
obtained by letting $\psi_2:=G$. As a first example, take $\d F(\tau)=\d
\tau$, $\varphi_{1}(\tau)=\tau^{\lambda-1}$, that is completely monotonic
for $\lambda \in (0,1)$, $g(\tau; \alpha_i, \nu_i)=
(1+\alpha(\ss_i)\tau)^{-\nu(\ss_i)}$, where $\alpha$ and $\nu$ are
supposed to be strictly positive functions of $\ss_i$ ($i=1,2$) and
additionally $0 < \alpha(\ss_i),\nu(\ss_i)< \pi$. One can readily verify
that all integrability conditions in Theorem \ref{nonstationary} are
satisfied. Letting
$k=\frac{|\Sigma(\ss_1)|^{1/4}|\Sigma(\ss_2)|^{1/4}}{|\Sigma(\ss_1,\ss_2)|^{1/2}}
Q (\ss_1,\ss_2)^{\lambda-1}$ and using $[3.259.3]$ in Gradshteyn and
Ryzhik (1980), one obtains the following covariance function
$${\rm C}_{s}(\ss_1,\ss_2)= k \alpha(\ss_1)^{-\lambda} {\rm
B}(\lambda, \nu(\ss_1)+\nu(\ss_2)-\lambda) _{2}F_{1}\left (\nu(\ss_2),
\lambda; \nu(\ss_1)+\nu(\ss_2); 1-\frac{\alpha(\ss_2)}{\alpha(\ss_1)}
\right ),
$$
where ${\rm B}(.,.)$ is the Beta function, and $_2F_1(.,.,.,.)$ is the
Gauss hypergeometric function.


Another example is obtained by considering $F(\d \tau)=\exp(-\tau)\d
\tau$, $\varphi_1(\tau)=\tau^{\nu-1}$, with $\nu \in (0,1)$, $g(\tau;
\ss_i)= \exp(-\frac{\alpha(\ss_i)}{2} \tau)$, with $\alpha(\ss_i)$
strictly positive ($i=1,2$) and using $[3.478.4]$ of Gradshteyn and Ryzhik
(1980), we find that
$${\rm C}_{s}(\ss_1,\ss_2)=2k \left (\frac{\alpha(\ss_1)+\alpha(\ss_2)}{2} \right )^{-\nu/2} \mathcal{K}_{\nu}\left (2 \left (\frac{\alpha(\ss_1)+\alpha(\ss_2)}{2}\right )^{1/2}\right ) $$
is a nonstationary spatial covariance that allows for local geometric
anisotropy, but has a fixed smoothing parameter.

It should be stressed that numerical integration under the setting
(\ref{nonstationarity_theorem}) could outperform previously proposed
models if the objective is to find a different type of interaction between
the local parameters characterising the integrating function $g$. In all
the examples proposed by Paciorek and Schervish (2006) and Stein (2005b)
the varying-smoothing, spatially adaptive parameter is obtained as the
semi-sum of a parameter acting on the location $\ss_1$ with another one
depending on the location $\ss_2$. This is a serious limitation of the
method, as only one type of interaction can be achieved. Starting from a
different setting, Pintore and Holmes (2004) obtain the same type of
spatially adaptive smoothing parameters. Thus, quasi-arithmetic
functionals could be of help, at least through numerical integration.
Finding a closed form for a functional different than $\mathcal{Q}_{G}$
remains an open problem to which we did not find any solution for the
meantime.

\section{Conclusions and discussion}

In this work, novel results are presented concerning permissible
spatial-temporal covariance functions  in terms of the theory of quasi
arithmetic means, and valuable insight is gained about their space-time
structure.  The theory of quasi arithmetic means is used, together with
the relevant permissibility criteria, to derive new classes of
nonseparable space-time covariances and to investigate their properties in
considerable detail.

There are several possible avenues for research based on the results of
the present paper.  From the
 spatial and spatiotemporal statistics perspectives, the QARF representation considered in this paper seems very promising.  E.g., the QARF may provide a
naive separability testing procedure, as follows.  Consider the set
$\Psi_{\rm cm}$ of all possible generators of the QARF class and let
$\Theta$ be the set of parameter vectors indexing the generators, {\em
i.e.} $\Theta: \{\boldsymbol{\theta} \in \RR^p :
\varphi_{\boldsymbol{\theta}} \in \Phi_{\rm cm}\}$. Thus, testing for
separability of the covariance function associated to the QARF is
equivalent to testing for the null hypothesis
$\varphi_{0}:=\varphi_{0}(x;\boldsymbol{\theta}_0)=\exp(-\boldsymbol{\theta}_{0}
x)$, $x>0$. This is a topic worth of further investigation.

Another interesting topic is that of generator estimation. A tempting
choice would be to consider the approach proposed by Genest and Rivest
(1993), and that in Ferguson {\em et al.} (2000) who find a serial version
of the Kendall's tau and a relationship between this concordance index and
the generating function $\varphi$. In this case, one would extend the
result of these authors at least to the lattice $\mathbb{Z}^{2}$.

Future research effort would also focus on a deeper study of the
properties of quasi-arithmetic covariances along several directions. For
instance, it would be interesting to find limit properties of the
covariance functions specified through this construction.

\section*{Acknowledgements} The authors are indebted to Christian Genest
for his unconditional support and useful discussions during the
preparation of this manuscript.

\section{Appendix}
\subsection{Proof of Proposition \ref{prop:archim-ordering}}

Let us begin by showing the final part, {\em i.e.} that
\begin{equation}
\mathcal{Q}_{\psi_1}(f_1,\ldots,f_n) \leq
\mathcal{A}_{\psi_2}(f_1,\ldots,f_n), \label{prova1}
\end{equation}
for $\psi_1,\psi_2 \in \mathfrak{Q}_{\rm cm}$, if and only if
$\varphi_{1}^{-1} \circ \varphi_2$ is subadditive.

Let $g = \varphi_{1}^{-1} \circ \varphi_{2}$ and denote
$x_i=\varphi_{2}^{-1} \circ f_i(\xx_i)$.

Thus, (\ref{prova1}) is equivalent to \begin{equation} \varphi_1 \left (
\sum_{i=1}^{n}\theta_i g(\xx_i) \right ) \leq \varphi_2 \left
(\sum_{i=1}^n \theta_i \xx_i \right ). \label{orderings_4}
\end{equation}
For the necessity, assume $\mathcal{A}_{\psi_1} \leq
\mathcal{A}_{\psi_2}$. Applying $\varphi_{1}^{-1}$ in both sides of
(\ref{orderings_4}) we obtain the result.

For the if part, assume (\ref{orderings_4}) holds. Thus, applying
$\varphi_1$ to both sides of the same inequality the result holds.

Let us now prove that $\mathcal{A}_{\Pi}(f_1,\ldots,f_n) \leq
\mathcal{A}_{\psi}(f_1, \ldots, f_n) \leq
    \sum \theta_i f_i$.

The lower bound in the inequality is direct consequence of Lemma 4.4.3,
Corollary 4.6.3 in Nelsen (1999, p. 110) and references therein, while the
upper bound is a direct property of convex functions.

\subsection{Proof of Proposition 2}
Recall that, for Bernstein's theorem (Feller, 1966, p.439), $\varphi \in
\Phi_{{\rm cm}}$ if and only if
\begin{equation}
\varphi(t)=\int_{0}^{\infty}e^{-r t} \d F(r) \label{bernstein},
\end{equation}
with $F$ a positive and bounded measure. Now, for the proof of {\bf (a)},
notice that, if $\varphi \in \Phi_{{\rm cm}}$, equation
(\ref{eq:christakos1}) can be written as $$\int_{0}^{\infty} \exp \left
(-r \sum_{i=1}^{n} \theta_i \varphi^{-1} \left ( {\rm C}_{i}(\hh_i) \right
) \right ) \d F(r).
$$ Now, observe that if, for every $i$, $\varphi^{-1} \circ {\rm C}_i$ is
a variogram, then so is the sum $\sum_{i=1}^{n} \theta_i \varphi^{-1}
\circ {\rm C}_{i}$. Then, by Schoenberg theorem (1939), we have that for
every positive $r$, the integrand in the formula above is a covariance
function. So is the positive scale mixture of covariances. This completes
the proof.

For {\bf (b)}, it is sufficient to notice that, by Schoenberg theorem
(1938), the mapping $\gamma : \RR^{d} \to \RR_{+}$ is the variogram
associated to an intrinsically stationary and isotropic random field if
and only if $\gamma(\hh)=\psi(\|\hh \|^2)$, for $\psi$ a Bernstein
function. The rest of the proof comes from the same arguments of point
{\bf (a)}.

For {\bf (c)}, notice that, being continuous,  increasing and concave on
$[0, \infty)$, each of the functions $\varphi^{-1} ({\rm C}_i)$, is
negative definite on $\RR$, according to a P{\' o}lya type criterion (see
Berg and Forst (1975), Proposition 10.6). So is their sum. Since $\varphi
\in \Phi_{{\rm cm}}$, by Schoenberg's theorem (cf. Berg and Forst, 1975),
we get once again the result.

\subsection{Proof of Proposition \ref{christakos4}}

We shall only prove the result in (\ref{eq:associativity0}), as
(\ref{eq:associativity}) and (\ref{eq:associativity2}) follow the same
argument. Recall that here we impose $\theta_i=1/n$ $\forall i$. Call
$\hh^{(1)}=(\| \hh_1 \|,\ldots,\| \hh_k \|)^{\prime}$ and $\hh^{(2)}=(\|
\hh_{k+1}\| ,\ldots,\| \hh_n\| )^{\prime}$, $k<n$. Also, let $f_1 = 1/n
\sum_{j=k+1}^{n} \varphi_1^{-1} \circ {\rm C_j}$, $f_2= \varphi_1^{-1}
\circ \varphi_2$ and $f_3= 1/n \sum_{i=1}^{k} \varphi_2^{-1} \circ {\rm
C}_i$. Thus, equation (\ref{eq:associativity0}) can be written, using
Bernstein's theorem, as $$\int_{0}^{\infty} \ee^{-r f_1(\hh^{(2)})-r f_2
\circ f_3 (\hh^{(1)}) } \d F_1(r)
$$
where $F_1$ is the distribution associated to its Laplace transform
$\varphi_1$. Thus, the proof follows straight by arguments of the previous
proofs.

\subsection{Proof of Theorem \ref{nonstationary}}
The result is a direct consequence of the work of Stein (2005b).
First, observe that ${\rm C}(\ss_1,\ss_1)=\int_{0}^{\infty} g(\tau; \ss_1)
\d F(\ss_1) < \infty$. Now, we need to show that $\sum_{i=1}^{n} a_i a_j
C(\ss_i,\ss_j) \geq 0$, for every finite system of arbitrary real
constants $a_i$, $i=1,\ldots,n$. Write $K_{i}^{\tau,r,r^{\prime}}$ for the
normal density centered at $\ss_i$ and with covariance matrix $(\tau r
r^{\prime})^{-1} \Sigma_{\ss_i}$. Also, by abuse of notation,
$\varphi_1:=\mathcal{L}[F_1]$ and $\varphi_2:=\mathcal{L}[F_2]$. Then, by
a convolution argument in Paciorek (2003, p.27), by Bernstein's theorem,
and using Fubini's theorem, we get
 {\tiny
\begin{eqnarray} &&\sum_{i,j=1}^{n} a_i a_j C(\ss_i,\ss_j)= \nonumber \\ &=& \pi^{p/2}\sum_{i,j=1}^{n} a_i a_j |\Sigma_{\ss_i}|^{\frac{1}{4}}|\Sigma_{\ss_j}|^{\frac{1}{4}}\int_{0}^{\infty} \int_{0}^{\infty}
\int_{0}^{\infty} \left (
\int_{\RR^p}K_{i}^{\tau,r,r^{\prime}}(\mathbf{u})
K_{j}^{\tau,r,r^{\prime}}(\mathbf{u}) \d \mathbf{u} \right )(\tau r
r^{\prime})^{-p/2} \ee^{-r^{\prime}(1/2\varphi_2^{-1}\circ g(\tau;
\ss_i)+1/2 \varphi_2^{-1} \circ g(\tau; \ss_j) )}\d F(\tau) \d F_1 (r) \d
F_2 (r^{\prime}) \nonumber \\
&=&\pi^{p/2} \int_{0}^{\infty} \int_{0}^{\infty} \int_{0}^{\infty}
\int_{\RR^p} \left ( a_i |\Sigma_{\ss_i}|^{\frac{1}{4}} K_i^{\tau, r,
r^{\prime}}(\mathbf{u}) \ee^{-r^{\prime} 1/2 \varphi_{2}^{-1} \circ
g(\tau;\ss_i)} \right )^2 \d \mathbf{u} (\tau r r^{\prime})^{- p / 2} \d
F(\tau) \d F_1(r) \d F_2 (r^{\prime}) \geq 0
\end{eqnarray}
}Thus, the proof is completed.

\subsection{M.S. differentiability of the QARF}

Recall that the existence of the $k$-th order $i$-th mean square partial
derivative of $Z$ is related to the existence of the $2k$-th order $i$-th
mean square partial derivative of the covariance function. This is in turn
related to the spectral moments through the following formula (Adler,
1981, p.31):

$$(-1)^k \frac{\delta^{2k} {\rm C}(\hh)}{\delta h_i^{2k}}|_{\hh=0} = \int_{\RR^d} \omega_i^{2k} \d \hat{{\rm C}}(\boldsymbol{\omega})<\infty. $$

Under the setting in Section 4.2, we ensure the existence of the spectral
density $\hat{\rm c}$, thus

\begin{eqnarray}
\int_{\RR^d} \omega_i^{2k} \hat{\rm c}(\boldsymbol{\omega}) \d
\boldsymbol{\omega} &\propto& \int_{\RR^d} \omega_i^{2k}\int_{\RR^d}
\ee^{-i \boldsymbol{\omega}^{\prime}
\hh} {\rm C}(\hh) \d \hh \d \boldsymbol{\omega} \nonumber \\
&=& \int_{\RR^d} \omega_i^{2k}\int_{\RR^d} \ee^{-i
\boldsymbol{\omega}^{\prime} \hh} \int_{0}^{\infty} \ee^{-r \sum_i
\theta_i \nu_i (|h_i|)} \d F(r) \d \hh \d
\boldsymbol{\omega} \nonumber \\
&=&\int_{\RR^d} \int_{0}^{\infty} \omega_i^{2k} \left (
\int_{\RR^{d-1}}\ee^{-i \tilde{\boldsymbol{\omega}}^{\prime} \tilde{\hh}}
\ee^{-r
 \sum_{j \neq i} \theta_j \nu_j (|h_j|)} \d \tilde{\hh} \right ) \left (
\int_{\RR} \ee^{-r \theta_i \nu_i (|h_i|)} \d h_i \right ) \d F(r) \d
\boldsymbol{\omega} \nonumber
\end{eqnarray}
where $\tilde{\hh}=(h_1, \ldots, h_{i-1}, h_{i+1}, \ldots, h_d)^{\prime}
\in \RR^{d-1}$ and where we repeatedly make use of Fubini's theorem for
standard integrability criteria. Now, observe that, for the assumption of
absolute integrability of $\exp(-r \theta_i \nu_i(x))$ for any positive
$r$, the last equality can be written as

$$  \int_{0}^{\infty} \int_{\RR^d} \omega_i^{2k} \hat{\rm c}(\tilde{\boldsymbol{\omega}};r) \hat{\rm c}(\omega_i;r) \d \boldsymbol{\omega} \d F(r)$$
where $\hat{{\rm c}}(\tilde{\boldsymbol{\omega}};r):=\mathcal{F}^{-1}[{\rm
C}_{r}](\tilde{\boldsymbol{\omega}})$, with ${\rm
C}_r(\tilde{\boldsymbol{h}})=\exp(-r \sum_{j \neq i} \theta_j \nu_j
(|h_j|))$ and where $\hat{\rm c}(\omega_i;r)=\mathcal{F}^{-1}[{\rm
C}_{r}](\omega_i)$, with ${\rm C}_{r}(h_i)=\exp(-r \theta_i \nu_i
(|h_i|))$. Then, by noticing that $\int_{\RR^{d-1}}\hat{\rm
c}(\tilde{\boldsymbol{\omega}};r) \d \tilde{\boldsymbol{\omega}}={\rm
C}_r(\tilde{\boldsymbol{0}}) < \infty$, one
obtains that the integral above is finite if and only if the function \\
$r \mapsto \int_{\RR} \omega_i^{2k} \hat{\rm c}(\omega_i;r) \d \omega_i$
is $F$-measurable. Thus, the proof is completed.

\end{document}